\magnification=\magstep1   
\hsize=15truecm            
\vsize=20truecm            
\parindent=2em             
\vskip2\baselineskip       
\catcode`\@=11             
\font\eightrm=cmr8         
\font\eightbf=cmbx8        
\newfam\msbfam             
\newfam\euffam             
\font\tenmsb=msbm10        
\font\sevenmsb=msbm7       
\font\teneuf=eufm10        
\font\seveneuf=eufm7       
\textfont\msbfam=\tenmsb   
\scriptfont\msbfam=\sevenmsb%
\textfont\euffam=\teneuf   
\scriptfont\euffam=
\seveneuf                  
\def\H{{\cal H}}    
\def\Q{{\fam\msbfam Q}}    
\def\N{{\fam\msbfam N}}    
\def\R{{\fam\msbfam R}}    
\def\Z{{\fam\msbfam Z}}    

\catcode`\@=12             
\def\.{{\cdot}}            
\def\<{\langle}            
\def\>{\rangle}            
\def\({\big(}              
\def\){\big)}              
\def\defi{\buildrel\rm def 
\over=}                    
\def\hat{\widehat}         
\def\implies{
\hbox{$\Rightarrow$}}      
\def\Homeo{{\cal H}}       
\def\id{\mathop{\rm id}    
\nolimits}                 
\def\im{\mathop{\rm im}    
\nolimits}                 
\def\pr{\mathop{\rm pr}    
\nolimits}                 
\def\Aut{\mathop{\rm Aut}  
\nolimits}                 
\def\defi{\buildrel\rm def 
\over=}                    
\def\ssk{\smallskip}       
\def\msk{\medskip}         
\def\bsk{\bigskip}         
\def\nin{\noindent}        
\def\cen{\centerline}      
\font\smc=cmcsc10          
\def\phi{\varphi}          
\newcount\litter                                         
\def\newlitter#1.{\advance\litter by 1                   
\edef#1{\number\litter}}                                 
\def\qedbox{\hbox{$\rlap{$\sqcap$}\sqcup$}}              
\def\qed{\nobreak\hfill\penalty250 \hbox{}               
\nobreak\hfill\qedbox\vskip1\baselineskip\rm}            
\def\arr{\hbox to 20pt{\rightarrowfill}}                 
\def\larr{\hbox to 20pt{\leftarrowfill}}                 
\def\mapdown#1{\Big\downarrow\rlap{$\vcenter{\hbox{$\scriptstyle#1$}}$}} 
\def\lmapdown#1{\llap{$\vcenter{\hbox{$\scriptstyle#1$}}$}\Big\downarrow}
\def\mapright#1{\smash{\mathop{\arr}\limits^{#1}}}
\def\lmapright#1{\smash{\mathop{\arr}\limits_{#1}}}
\def\mapleft#1{\smash{\mathop{\larr}\limits^{#1}}}
\def\lmapleft#1{\smash{\mathop{\larr}\limits_{#1}}}

\long\def\alert#1{\parindent2em\smallskip\line{\hskip\parindent\vrule%
\vbox{\advance\hsize-2\parindent\hrule\smallskip\parindent.4\parindent%
\narrower\noindent#1\smallskip\hrule}\vrule\hfill}\smallskip\parindent0pt}

\newlitter\andii.
\newlitter\arens.
\newlitter\arh.
\newlitter\arch.
\newlitter\bourb.
\newlitter\bred.
\newlitter\cook.
\newlitter\gart.
\newlitter\grooti.
\newlitter\grootii.
\newlitter\grootiii.

\newlitter\dijk.
\newlitter\droste.
\newlitter\ellis.
\newlitter\hofi.

\newlitter\amt.

\newlitter\keesling.
\newlitter\keesii.

\newlitter\ryb.

\newlitter\ttt.

\def\compbook{\hofi}

\def\epsilon{\varepsilon}

\cen{\bf Representing  a Profinite Group}

\cen{\bf as the Homeomorphism Group of a Continuum }
\msk
\cen{ by Karl H. Hofmann and Sidney A. Morris}
\bsk\noindent
{\eightbf Abstract}. {\eightrm 
We contribute some information towards finding a general algorithm for constructing, for a given profinite group $\scriptstyle G$, a compact connected space $X$ such that the full homeomorphism group ${\scriptstyle \cal H}\scriptstyle(X)$ with the compact-open topology is isomorphic to $\scriptstyle G$ as a topological group. 
It is proposed that one should find a {\it compact topological} oriented
graph 
$\scriptstyle \Gamma$ such that 
$\scriptstyle G\cong\Aut(\Gamma)$. 
The replacement of the edges of 
$\scriptstyle \Gamma$ by rigid continua should work as is exemplified in
various instances where discrete graphs were used.

It is shown here that the strategy can be implemented for 
profinite monothetic groups 
$\scriptstyle G$.


}

\bsk

\cen\nin{\smc 1. Introduction}
\msk
\noindent
One knows that the compact homeomorphism group $\H(X)$ of a Tychonoff space
has to be profinite ([\amt], [\keesling]). In the converse direction
{\smc Gartside} and  {\smc Glyn} [\gart] have established that every metric
profinite group is the homeomorphism group of a continuum (i.e.\ a compact connected metric space).

For the goal of representing a given group as the homeomorphism group of
a space, authors have pursued the following strategy:

Step (1): find some connected graph $\Gamma$, usually oriented, and find
an isomorphic representation $\pi\colon G\to \Aut(\Gamma)$; the
standard attempt is to use some form of  Cayley graphs (see [\grooti],
[\gart], [\arch])

Step (2): find a rigid continuum $C$, that is, a continuum, that is,
compact connected metric space,  
whose only continuous selfmaps are the identity and the constant 
function (see [\cook], [\grootiii])
and replace each of the directed edges of $\Gamma$ by $C$ or a variant
obtaining a connected space $X$; finally obtain an
isomorphism  $\gamma\colon\Aut(\Gamma) \to \H(X)$
(see [\grooti], [\grootiii], [\arch]). Obtain an isomorphism
$\gamma\circ\pi\colon G\to \H(X)$.

\msk 
All known variations of the strategy are highly technical, and 
different variations lead to rather different phase spaces $X$.
It would be nice to find a construction which is in some way
canonical, perhaps even functorial. However,
one of the major obstructions for
representations of a profinite group in a combination
with  graph theoretical methods is
that homeomorphism groups, like all automorphism groups in a category
are, in no visible way, functorial.
 
\bsk

We propose, that in Step (1) one should in fact go more than halfway
and construct a compact connected directed graph $\Gamma$ 
and {\it then} apply Step (2) to achieve the final goal.

\msk
In the following we show that such constructions are possible in
principle and yield for every profinite monothetic group $G$ 
a continuum $X$ such that $\H(X) \cong G$. while not all compact 
monothetic groups are metric, the profinite ones among them are.

Thus, in the vein of a general existence result, our construction
yields nothing new beyond what {\smc Gartside} and {\smc Glyn}
have shown in [\gart]. However, the continua we construct are completely 
different from those produced in [\gart] and the proposed construction may
turn out to be useful in the future.

\bsk


\bsk

\cen{\smc 2.   Directed topological graphs.}
\msk
In order to construct topological spaces with prescribed homeomorphism
groups we first construct directed topological graphs 
with prescribed automorphism
groups.
\msk
\bf\nin 
Definition 2.1. \quad \rm
A  {\it directed (topological) graph} 
is a triple $\Gamma=(V,E,\eta)$ consisting of topological spaces
 $V$ and $E$ and a continuous function
$\eta\colon E\to V\times V$ such that
$$V=\im(\pr_1\circ \eta)\cup\im(\pr_2\circ \eta). \leqno(\dag)$$  

The set $V$ is called the space of {\it vertices} and $E$ is called
the space of {\it (oriented) edges}. If $e\in E$  we write
$\eta(e)=(e^1,e^2)\in V\times V$, 
then $e^1$ is the {\it origin}
and $e^2$ is the {\it target} of $e$. Condition $(\dag)$ says
$$ V=\bigcup_{e\in E}\{e_1,e_2\},\leqno{(\ddag)}$$
and this means that there is no vertex that is not an endpoint of an
edge.  Note that 
we allow a whole space $\eta^{-1}(v_1,v_2)$
of (directed) edges from $v_1$ to $v_2$.
We shall, however, not use this fact in the sequel.

If the spaces $V$ and $E$ of a directed graph $\Gamma=(V,E,\eta)$
are discrete, we recover the more classical concept of a directed graph. 

\msk

\bf\nin
Example 2.2.\quad\rm (i) Let $n$ be natural number $n>2$ and let
$\Z(n)=\Z/n\.\Z$ be the cyclic group of $n$ elements. Define
$$\eqalign{%
V&=\Z(n)=\{m+\Z\in \Z/n\Z: m=0,1,\dots,n-1\},\cr
E&=V,\cr
\eta(m+n\.\Z)&=(m+n\.\Z, m+1+n\.\Z)\in V\times V.\cr}$$

The directed graph $C(\Z(n))\defi (V,E,\eta)$ is the $n$-{\it cycle}. 
It is the Cayley-graph of the pair
$(\Z(n), \{1+n\Z\})$ consisting of the cyclic group of order
$n$ and the singleton generating set containing the element
$1+n\Z$.

\msk (ii) More generally, let $G$ be a topological  group
and  let $g\in G$. We set
$$\eqalign{%
V&=G,\cr
E&=V,\cr
\eta(x)&=(x, xg)\in V\times V.\cr}$$

The directed graph $C(G)\defi(V,E,\eta)$ is the {\it topological
Cayley graph} of $(G,\{g\})$.   

\msk
(iii) Taking  $\Z$ with the discrete topology 
we obtain the {\it Cayley graph} $C(\Z)$  of $(\Z,\{1\})$,
the chain $\Z$ with its natural order-orientation. \qed

\msk
\bf\nin
Definition 2.3. \quad \rm A doubly pointed 
connected topological space 
$${\bf L}=(L,b_1,b_2) \hbox{ with }b_1\ne b_2$$
is called a {\it link}.\qed

Typically {\bf I}$=([0,1],0,1)$ is a link:
an interval joining its endpoints. If {\bf L}$= (L,b_1,b_2)$ is
a link and a Tychonoff space, then there is a continuous function
$F\colon{\bf L}\to {\bf I}$, $F\colon L\to I$, $F(b_1)=0$ and
$F(b_2)=1$, called a {\it morphism of links}.
\msk
\bf\nin
Construction 2.4.\quad \rm Let $\Gamma=(V,E,\eta)$ be a topological
directed graph and {\bf L}$=(L,b_1,b_2)$ a link.
We construct a
 space $|\Gamma|_{\bf L}$ from these data by ``inserting into
each oriented edge $e=\eta(v)=(e^1,e^2)$
 a copy of $\{e\}\times L$ of the link  {\bf L} such that 
$(e,b_1)$ is identified with $e^1$ and $v$ 
while  $(e^2,b_2)$ is identified with with $e^2$ and $v$.''

Indeed we let $X=E\times L$ 
and define an equivalence
relation $\rho$ on $X$ with the following equivalence classes:
$$
\rho(e,x) =\cases{\{(e,x)\} &if $x\ne b_1,b_2$,\cr 
          \{(e',b_1):\eta(e')=(v,(e')^2)\}\cup \{v\} &if 
                                          $\eta(e)=(v,e^2),\ x=b_1$, \cr
           \{(e'',b_2):\eta(e'')= ((e'')^1,v)\}\cup \{v\} &if  
                                          $\eta(e)=(e^1,v),\  x=b_2$, \cr
          \{(e^*,x):\eta(e^*)=(v,v)\}\cup \{v\}&if $\eta(e)=(v,v),\  
                                              x\in\{b_1,b_2\}$.\cr}
$$

Then let $|\Gamma|_{\bf L}\defi X/\rho$.  The space 
$|\Gamma|_{\bf L}$ is called the 
{\it topological realisation} 
of $\Gamma$ via the link {\bf L}$=(L,b_1,b_2)$.

If {\bf L} is a Tychonoff link and $F\colon{\bf L}\to{\bf I}$
a morphism of links then our construction obviously induces
a morphism $F^*\colon \Gamma_{\bf L}\to\Gamma_{\bf I}$ of 
topological realisations.\qed

 Notice that in the case of a Cayley graph of a 
group $G$ with an element $g\in G$,
 the quotient space $((E\times L))/\rho$ can
be expressed in the form 
$$|C(G)|_{\bf L}=(G\times L)/\rho.\leqno(*)$$ and that there
is a morphism $F^*\colon |C(G)|_{\bf L}\to |C(G)|_{\bf I}$ of
topological realisations given by $F^*(\rho_L(a,x))=\rho_I(a,f(x))$.
\msk

The verification of the details of the following examples
is straightforward.
\msk

\bf\nin
 Example 2.5.\quad\rm  $|C(\Z(n))|_{\bf I}$ is a circle and 
$|C(\Z)|_{\bf I}$ 
is homeomorphic to $\R$.\qed

\msk
For the following, let
 $R$ be a  compact connected space
and pick two different points $b_1, b_2\in E$.
and {\bf R}$=(R,b_1,b_2)$ the corresponding link.
\msk
\bf\nin
Lemma 2.6. \quad\it Let $(X,x_0)$ be a compact connected pointed space and
 \hfill \break $(R,b_1,b_2)$  a doubly pointed  de Groot-continuum. 
Assume that $X$ and $R$ are 
disjoint  with the exception of $b_2$ and $x_0$
which are assumed to be equal. Then a continuous function 
$f\colon R\to X\cup R$ is exactly one of the following kind

{\parindent3em

\item{\rm(i)} $f$ is constant.

\item{\rm(ii)} $f(R)=R$ and its corestriction $R\to f(R)$ is the
identity map of $R$.

\item{\rm(iii)} $f(R)\subseteq X$ and $f(b_2)=b_2$.

\item{\rm(iv)} $f(R) \cap R=\emptyset$. 

}\rm

\msk\nin
\bf Proof. \quad\rm
The function $\pi\colon R\cup X\to R$ defined by 
$\pi(X)=\{b_2\}$ is continuous; hence the continuous self-map
$\pi\circ f\colon R\to R$ is either the identity or is constant
with image $\{b_2\}$.\qed

\msk
\bf\nin
Example 2.7.\quad\rm 
Let ${\bf R}=(R,b_1,b_2)$ be a doubly pointed de Groot-continuum and
$F\colon {\bf R}\to{\bf I}$ a morphisn of links. 
 Then

(i) $|C(\Z(n))|_{\bf R}$ is a one-dimensional continuum $X$ whose homeomorphism
group $\Homeo(X)$ is isomorphic to the cyclic group $\Z/n\.\Z$.
Moreover $F^*\colon|C(\Z(n))|_{\bf L}\to|C(\Z(n))|_{\bf I}$ is a morphism
of realisations onto, or ``over'', a circle.
\ssk 

(ii)  $|C(\Z)|_{\bf R}$    
is a one-dimensional connected locally compact space
$X$ whose homeomorphism space $\Homeo(X)$ is isomorphic to the 
infinite cyclic group $\Z$. Moreover 
$F^*\colon |C(\Z)|_{\bf L}\to|C(\Z)|_{\bf I}$ is a morphism
of realisations over a line.
\msk
\bf\nin
Proof. \quad \rm From $(*)$ above recall that for a cyclic
group $Z=\Z/m\Z$, $n=0,1,\dots$ we have
$|C(Z)|_{\bf R}=(Z\times R)/\rho$ and  note via Lemma 2.6,
that  the action 
 $(n,(k+m\Z,x))\mapsto (n+k+m\Z,x):\Z\times Z\times R\to Z\times R$ 
gives a unique
action 
$(n, \rho(k+m\Z,x))\mapsto \rho(n+k+m\Z,x):
\Z\times |C(Z)|_{\bf R}\to |C(Z)|_{\bf R}$
representing the action of $\Homeo(|C(Z)|_{\bf R})$ on $|C(Z)|_{\bf R}$. \qed

We now generalize  Examples 2.7(i) to monothetic compact groups by
utilizing Example 2.7(ii).

\msk\bf\nin
Main Lemma 2.8. \quad\rm
Let $G$ be a compact 
group with a nonidentity element $g$. 
Let $C(G)$ be the topological Cayley graph of $(G,g)$. Let
{\bf L}$=(L,b_1,b_2)$ be  a compact link. Let $F\colon{\bf L}\to{\bf I}$
be a morphism of links onto the interval.
Then the following conclusions hold:

{\parindent2.5em

\item{\rm(i)} $\Z$ acts freely with discrete orbits 
on $G\times |C(\Z)|_{\bf L}$ via $m\.(a,x)=(a^m, -m\.x)$ and  
the compact orbit space $(G\times |C(\Z)|_{\bf L})/\Z$
is naturally isomorphic to  $|C(G)|_{\bf L}$.
It is locally homeomorphic to
$G\times |C(\Z)|_{\bf L}$ under the orbit map of the $\Z$-action.

\item{\rm(ii)} If $G$ is monothetic with generator $g$, then
$|C(G)|_{\bf L}$ is a compact connected 
 Hausdorff space. If $G$ is profinite monothetic, then
$|C(G)|_{\bf I}$ is homeomorphic to a solenoid ($p$-adic if
$G=\Z_p$, the additive group of $p$-adic integers), and
$F^*\colon |C(G)|_{\bf L}\to |C(G)|_{\bf I}$ is a morphism of
realisations over a solenoid.

\item{\rm(iii)} The group $G$ acts on $G\times |C(\Z)|_{\bf L}$
by the left regular action on the left factor. This 
action commutes with the $\Z$-action.
It thus induces an action 
$$(a, \Z\.(b,x))\mapsto \Z\.(a+b,x): G\times (G\times |C(\Z)|_{\bf L})/\Z
\to G\times |C(\Z)|_{\bf L})/\Z.$$
The orbit space $|C(G)|_{\bf L}/G$ is homeomorphic to
the space $L/\{b_1,b_2\}$ arising from $L$ by
identifying the two points $\{b_1,b_2\}$.

}
\msk
\bf\nin
Proof. \quad \rm
(i) The assertions on the action are 
straightforward.
Now $$|C(\Z)|_{\bf L}= (\Z\times L)/\rho$$
has $\dot L\defi\rho^{-1}(\{0\}\times L)/\rho\cong L$
as a fundamental domain for the $\Z$ action, i.e. 
each orbit meets $\dot L$ only once except the orbit
of $b_1$ which meets $\dot L$  in $\{\dot b_1,\dot b_2\}$.
Now we deduce that in this spirit  
$G \times \dot L$ is a fundamental domain
of the $\Z$-action on $G\times|C(\Z)|_{\bf L}$ 
Thus the orbit space $(G\times |C(\Z)|_{\bf L})/\Z$ 
is a continuous image of 
$G\times \dot L$ and therefore is compact.
\ssk
Next we denote by $\sigma$ the equivalence relation
on $\Z\times L$ which identifies
$(n,b_2)$ and $(n+1,b_1)$ for all $n\in\Z$ so that
$|C(\Z)|_{\bf L}=(\Z\times L)/\sigma$. We let
$q\colon G\times |C(G)|_{\bf L}\to(G\times |C(G)|_{\bf L})/\Z$
be the orbit map of the $\Z$ action. If 
$\Z\.(g,m,x)=\{(g^n,m-n,x):n\in\Z\}$ we set $\alpha(\Z\.(g,m,x))
=\{(g^n,\sigma≈(m-n,x)):n\in\Z$.
Further we let $\pi\colon (G\times\Z)/\Z\to G$ be the isomorphism
given by $\pi(g) =(g,z)\Delta$ for $\Delta=\{(g^n,-n): n\in\Z\}$.
Let $\rho$ on $G\times L$ be the equivalence relation which 
collapses $(ag,b_1)$ and $(a,b_2)$ for all $a\in G$. Also, we define
$\beta\colon{G\times\Z\over\Z}\times L\to{G\times L\over\Z}$ by
$\beta(\{(ag^n,m-n):n\in\Z\},x)=\beta(\pi(a),x)=\rho(a,x)$
and $\gamma\colon{G\times L\over\rho}\to{G\times|C(\Z)|_{\bf L}\over\Z}$
by $\gamma(\rho(a,x))=\{(ag^n,\sigma(m-n,x)):n\in\Z\}$. 

We then have a commutative diagram
$$\matrix{%
G\times \Z\times L&\mapright{Q}&{G\times \Z\times L\over\Z}&\mapleft{\Gamma}&
                       {G\times \Z\over \Z}\times L
                      &\mapleft{\pi\times\id_L}&G\times L\cr
\lmapdown{\id_G\times \sigma}&&\lmapdown{\alpha}&& 
              \mapdown{\beta}&&\mapdown{\rho}\cr
G\times |C(\Z)|_{\bf L}&\lmapright{q}&{G\times|C(\Z)|_{\bf L}\over\Z} &
                        \lmapleft{\gamma}&{G\times L\over\rho}&=&
                                            |C(G)|_{\bf L}\cr.}$$ 

Indeed the commuting of the first rectangle is clear from the definition
of the action of $\Z$ on $G\times|C(\Z)|_{\bf L}$, and the commuting
of the right rectangle is an immediate consequence of the definitions. 

The middle rectangle, however, commutes since for all $a\in G$ and $m\in\Z$
we have 
$$\leqalignno{%
&\alpha\circ \Gamma(\{(ag^n,m-n):n\in\Z\},b_2)=
                          \alpha(\{ag^n,m-n,b_2):n\in\N\})\cr
=&\{(ag^n,\sigma(m-n,b_2)):n\in\Z\}=\{(ag^n,\sigma(m-n+1,b_1)):n\in\Z\}&(1)\cr}$$
while
$$\leqalignno{%
&\gamma\circ \beta(\{(ag^n,m-n):n\in\Z\},b_2)=\gamma\rho(a,b_2)
=\gamma(\rho(ag,b_1))\cr 
=&\{(ag^{p+1},\sigma(m-p,b_1)):p\in\Z\}
=\{(ag^n,\sigma(m-(n-1),b_1)):n\in\Z\}.&(2)\cr}$$
Since (1) and (2) are obviously equal, the commuting of the middle
rectangle follows. We see at once that $\gamma$ is surjective since
$\alpha$ is surjective. 
We notice that $\gamma(\rho(a,x))=\gamma(\rho(a',x'))$ iff
$\{(ag^n,\sigma(m-n,x)):n\in\Z\}=\{(a'g^{n'},\sigma(m-n'),x'):n'\in\Z\}$
If $x\notin\{b_1,b_2\}$ these $\Z$=orbits on  agree iff and only
if their intersections with the fundamental domain 
$G\times \dot L$ agree. But $\sigma(k,x)\in \dot L$
iff one of the three possibilities apply: (i) $k=0$,
(ii)  $k=1$ and $x=b_1$, or (iii) $k=-1$ and $x=b_2$.
In the first case, $n=m=n'$ and $a=a'$, $x=x'$ follow.
In the second  case $n=m+1=n'+1$ is a possibility,
whence $ag^{m+1}=a'g^m$ that is $a'=ag$, $x'=b_1$ and $x=b_2$,
which implies that $\rho(a,x)=\rho(a',x')$. 
The other cases are discussed similarly and likewise yield
$\rho(a,x)=\rho(a',x')$. Therefore the continuous
function $\gamma$ between compact Hausdorff spaces is 
 bijective and therefore is a homeomorphism.
This completes the proof of (i)

\msk
(ii) First we have to prove connectivity of $|C(G)|_{\bf L}$. 
The space $ A=|C(\Z)|_{\bf L}$
is arcwise connected. Since $g^\Z$ is dense in $G$, the image 
$\Z\.|C(\Z)|_{\bf L}/\Z$ of $A$ is dense in 
$(G\times|C(\Z)|_{\bf L})/\Z$ which is naturally homeomorphic
$|C(G)|_{\bf L}$. Hence the latter space is connected.

Now let $A$ be a subgroup of the discrete group $\Q$
containing $\Z$ such that $T\defi A/\Z$ is infinite. Then
the character group $S\defi\hat A$ is a compact connected
one-dimensional abelian group called a {\it solenoid}.
The character group $G=\hat T$ is profinite monothetic.

By 2.7(ii) and (i)
we know that $|C(G)|_{\bf I}$ is homeomorophic to the
quotient group $(G\times\R)/\Delta$ for the subgroup
$\Delta=\{(n,-n):n\in\Z\}$. This quotient is one of the 
well-known representations of the solenoid $S$.
(See e.g.\ [\compbook],  Exercise E1.11., Theorem 8.22.)

\msk
(iii) We identify $|C(G)|_{\bf L}$ with $(G\times |C(\Z)|_{\bf L})/\Z$.  
The assertions are straightforward. 
The orbit space $|C(G)|_{\bf L}/G=
\((G\times|C(\Z)|_{\bf L})/\Z\)/G$ is isomorphic to
$\((G\times|C(\Z)|_{\bf L})/G\)/\Z\cong |C(\Z)|_{\bf L}/\Z\cong
L/\{b_1,b_2\}$. \qed

\msk\bf\nin
Theorem 2.9.\quad \it For any profinite monothetic group $G$
there is a compact connected 1-dimensional space $X$ such that
$\Homeo(X)\cong G$.

\msk\nin
\bf Proof.\quad \rm Let $g\in G$ be a generator of $G$. If $g$ has
finite order, there is nothing to prove because the assertion was established
in Example 2.7(i). We therefore assume for the rest of the proof 
that $g$ has infinite order. 
We apply Main Lemma 2.8 with $X=|C(G)|_{\bf R}$
for a doubly pointed de Groot continuum {\bf R}. 
By the Main Lemma 2.8 we know that $X$ is a compact connected 
space which is locally homeomorphic to the 
space $G\times |C(\Z)|_{\bf R}$ which is one-dimensional
since $R$ is one-dimensional.  Therefore $X$ is one-dimensional.
\msk 
We have to prove that $\Homeo(X)\cong G$.
This is the most delicate portion of the proof.
We identify $X$ with $(G\times|C(\Z)|_{\bf R})/\Z$.
For $h\in G$ let $\gamma_h\colon X\to X$ be defined by
$\gamma_h(\xi)=h\.\xi$ for the $G$-action on $P$.  
Let $\phi$ be a homeomorphism of $X$. We claim 
that $\phi$ is of the form $\gamma_h$ for some $h\in G$,
that is $\gamma_h(\Z\.(a,x))= \Z.(h+a,x)$.
 
The path components of $G\times |C(\Z)|_{\bf R}$ are
the spaces  $\{a\}\times |C(\Z)|_{\bf R}$ since 
$G$ is totally disconnected, and they are permuted by the 
action of $G$ on the left factor by the regular representation.
We may therefore consider $a=0$ without loss of generality.
The orbit 
$O\defi\Z\.(0,\sigma(0,b_1))=\{(n\.g,\sigma(-n,b_1)):n\in\Z\}$
meets $\{0\}\times |C(\Z)|_{\bf R}$ in an element
$(n\.g,\sigma (-n,b_1))$ only if $n\.g=0$ iff $n=0$
since $g$ has infinite order. It follows that the
orbit map 
$q\colon G\times |C(\Z)|_{\bf R}\to X$ maps each set
$\{a\}\times |C(\Z)|_{\bf R}$ continuously and bijectively 
onto 
$$Z_a\defi{\Z.(\{a\}\times |C(\Z)|_{\bf R})\over\Z}
={\Z\.a\times |C(\Z)|_{\bf R}\over\Z}$$
for $a\in G$. In particular, each of the sets $Z_a$ is
arcwise connected. 
We claim that the $Z_a$ are the arc components. This is
well-known in the  case of the solenoid $S=|C(G)|_{\bf I}
\cong (G\times \R)/\Delta\cong \hat A$ 
(cf.\ 2.8(ii) and its proof and [\compbook], Theorem 8.30).
We have the morphism $F^*\colon |C(G)|_{\bf R}\to|C(G)|_{\bf I}=S$
mapping $Z_0$ to $$(\{0\}\times \R)\Delta/\Delta=(\Z\times \R)/\Delta$$
the identity arc component of the solenoid $S$. Different 
sets $Z_a$ and $Z_b$ are mapped to different arc components in $S$
and so there can be no arc connecting a point in $Z_a$ to a point in
$Z_b$. Hence the $Z_a$ are precisely the arc components of 
$|C(G)|_{\bf R}$.

Thus the homeomorphism $\phi$ permutes the sets $Z_a$

We consider the particular arc component 
$$Z=Z_0=\Z.(\{0\}\times |C(\Z)|_{\bf R})/\Z.$$
Assume that $\phi(Z)=Z_a$. Then
$\gamma_a^{-1}\phi$ is a homeomorphism of $X$ which
maps the arc component $Z$ into itself.

The function $\epsilon\colon |C(\Z)|_{\bf R}\to Z$, 
$\epsilon(x)=\Z\.\sigma(0,x)$
is a continuous bijection.

We shall now invoke the arc component topology attached functorially
to a topological space as summarized in [\compbook], Appendix 4, p.\ 781.
Since $R$ and thus $|C(\Z)|_{\bf R}$ are locally arcwise connected,
the bijective function $\epsilon |C(\Z)|_{\bf R}\to Z$ is  the universal
map $\epsilon_Z\colon Z^\alpha \to Z$ of Lemma A4.1
of [\hofi], p.~781. Therefore, by Lemma A4.1(iv) every
homeomorphism of $Z$ lifts uniquely to a homeomorphism
of $|C(\Z)|_{\bf R}$ and thus is an action of $m\in\Z$ on $|C(\Z)|_{\bf R}$ by
Example 2.7(ii).

By the definition of the $\Z$-action on $G\times|C(\Z)|_{\bf R}$
according to 2.8(i), the action by $m$ on $|C(\Z)|_{\bf R}$ when
pushed down to $Z$ is induced  by the action of 
$$\gamma_m\colon X\to X,\quad
X=(G\times |C(\Z)|_{\bf R})/\Z.$$ 
Therefore the homeomorphism
$\gamma_{m-a}\phi=\gamma_m\gamma_a^{-1}\phi$ fixes the 
arc component $Z$ elementwise. However, $Z$ is dense in $X$.
Hence $\gamma_{m-a}\phi=\id_X$. Thus $\phi=\gamma_{a-m}$
and this had to be shown.\qed

As we have noted in the proof of 2.8(ii),
 a compact profinite group $G$ is monothetic, iff its
character group $\hat G$ is isomorphic to a subgroup  $A/\Z$
of the group $\Q/\Z$. The solenoid attached to this
monothetic group is the character group $\hat A$ of 
$A\subseteq \Q$.

\bsk

\bsk
Our feeling is that the occurrence of compact homeomorphism groups,
given certain restrictions, is not so rare as one might 
surmise initially even though the construction of compact 
spaces having a given homeomorphism group requires work. 
Ideally, one should be able to prove the following
\msk
\bf\nin
Conjecture.\quad \it Let $G$ be a compact group. Then the following
conditions are equivalent:
{\parindent2em

\item{\rm(1)} There is a compact  connected 
              space $X$ such that $G\cong\Homeo(X)$.

\item{\rm(2)} There is a compact space $X$ such that 
              $G\cong\Homeo(X)$.

\item{\rm(3)} $G$ is profinite.  

}
\rm
Note that we have (1)\implies(trivially) (2)\implies (3);
the implication (3)\implies (1) we have proved only for compact
monothetic groups $G$ and {\smc Gartside} and {\smc Glyn}
have proved it for arbitrary metric profinite groups. 
Cantor groups $\Z(2)^X$
(with arbitrary exponent $S$) Keesling [\keesii] was
able to represent as homeomorphism groups of one-dimensional
metric spaces.

\bsk

\nin{\smc 3. References} 

\msk\rm
{\parindent1.7em

\item{[\andii]}
Anderson, R. D.,
  The algebraic simplicity of certain groups of homeomorphisms,
    Amer. J. Math. {\bf80}  (1958),   955-–963.

\ssk

\item{[\arch]} 
Arhangel'skii, A. V, and M.~Tkachenko,
 Topological Groups and Related Structures,
  Atlantis Studies in Mathematics {\bf1}, 2008, 781pp. 

\ssk

\item{[\arens]}
Arens, R. F.,
  Topologies for homeomorphism groups,
    Amer. J. Math. {\bf68}   (1946),   593-–610.

\ssk

\item{[\bourb]}
Bourbaki, N., 
  Topologie g\'en\'erale,
    many publishers from Hermann,  Paris, ca 1950, to Springer Berlin, etc.,
     ca 2000.

\ssk

\item{[\bred]} 
Bredon, G., 
  Introduction to Compact Transformation Groups, 
    Academic \break Press, New York, 1972.

\ssk
\item{[\cook]} 
Cook, H., 
  Continua which admit only the identity mapping onto 
  non-degener\-ate subcontinua,
    Fund. Math. {\bf60} (1967), 241--249.

\ssk

\item{[\gart]}
Gartside, P. and A. Glyn,
   Autohomeomorphism groups,
     Topology Appl. {\bf129} (2003), 103--110.

\ssk

\item{[\grooti]}
de Groot, J.,
  Groups represented by homeomorphism groups.
    Math. Ann. {\bf138} (1959), 80--102. 

\ssk

\item{[\grootii]}
de Groot, J., and R. H. McDowell,
  Autohomeomorphism groups of $0$-dimens\-io\-nal spaces.
    Compositio Math. {\bf15} (1963), 203--209. 

\ssk

\item{[\grootiii]}
de Groot, J., and R. J. Wille, 
  Rigid continua and topological group-pictures.
    Archiv d. Math. {\bf9} (1958), 441--446.

\ssk

\item{[\dijk]}
Dijkstra, J. J.,  and J. van Mill, 
  On the group of homeomorphisms of the real 
  line that map the pseudoboundary onto itself, 
    Canad. J. Math. {\bf58} (2006), 529--547.

\ssk

\item{[\droste]} Droste, M., and R. G\"obel,
 On the Homeomorphism Groups of Cantor's Discontinuum and the Spaces
of Rational and Irrational Numbers, 
Bulletin of the London Mathematical Society {\bf34} (2002), 474--478.

\ssk

\item{[\ellis]} Ellis, R.,
 Locally compact transformation groups,
    Duke Math. J. {\bf24} (1957), 119--126

\ssk

\item{[\hofi]} 
Hofmann, K. H., and S. A. Morris,
  The Structure of Compact Groups, 
    Verlag Walter De  Gruyter Berlin, 1998, xvii+834pp.
    Second Revised and Augmented Edition 2006, xviii+858pp.

\ssk

\item{[\amt]}
---,
Compact Homeomorphism Groups are Profinite,
Preprint

\ssk

\item{[\keesling]} 
Keesling, J.,
   Locally compact full homeomorphism groups are zero dimensional,
     Proc. Amer. Math. Soc. {\bf29} (1971), 390--396.

\ssk

\item{[\keesii]}
---,
  The group of homeomorphisms of a solenoid,
Trans. Amer. Math. Soc. {\bf172} (1972), 390--396.

\ssk

\item{[\ryb]}
  Rybicky, T.,
   Commutators of homeomorphisms of a manifold,
     Universitatis Jagellonicae Acta Math {\bf23} (1996),
     153--1960.

\ssk
\item{[\ttt]}
tom Dieck, T., 
   Transformation Groups,
     Verlag Walter De  Gruyter Berlin, 1987, x+312pp.

}

\bsk\bsk
\line{%
\vbox{\hsize.48\hsize \baselineskip8pt \eightrm\obeylines
Karl H Hofmann
Fachbereich Mathematik
Technische Universit\"at Darmstadt
Schlossgartenstrasse 7
64289 Darmstadt, Germany
hofmann@mathematik.tu-darmstadt.de
\vglue3\baselineskip}\hfill
\vbox{\hsize.51\hsize \baselineskip8pt\eightrm\obeylines
Sidney A. Morris
School of Science, IT, and Engineering
University of Ballarat
P.O. 663, Ballarat
Victoria 3353, Australia, and
School of Engineering 
and Mathematical Sciences 
La Trobe University,
Bundoora Victoria 3086, Australia
morris.sidney@gmail.com}
}

\bye